\documentclass[11pt,a4paper]{amsart}
\usepackage{amsfonts,amscd,amssymb,amsmath,amsthm,hyperref,mathrsfs,xcolor,lscape,longtable}
\usepackage{microtype}
\newcommand{\bb}[1]{\mathbb{#1}}

\newcommand{\C}{\bb{C}}
\newcommand{\Q}{\bb{Q}}
\newcommand{\R}{\bb{R}}
\newcommand{\Z}{\bb{Z}}
\newcommand{\N}{\bb{N}}

\newtheorem{thm}{Theorem}[section]

\newtheorem{defn}[thm]{Definition}

\theoremstyle{definition}
\newtheorem{rmk}[thm]{Remark}

\theoremstyle{definition}

\numberwithin{equation}{section}
\newcommand{\Sp}{\mathrm{Sp}}

\newcommand{\GL}{\mathrm{GL}}

\begin{document}
\subjclass[2000]{??, ??}
\title{Thinness of some Hypergeometric Groups in $\Sp(6)$}
\author{Sandip Singh and Shashank Vikram Singh}
\address{Department of Mathematics, Indian Institute of Technology Bombay, Mumbai, India}
\email{sandip@math.iitb.ac.in}
\address{Department of Mathematics, Indian Institute of Technology Bombay, Mumbai, India}
\email{shashankvikrams6@gmail.com}
\begin{abstract}
We show the thinness of $7$ of the $40$ hypergeometric groups having a maximally unipotent monodromy in $\Sp(6)$. 
\end{abstract}
\maketitle
\section{Introduction}
A hypergeometric group $\Gamma(\alpha,\beta)$ associated to the pair of parameters $\alpha=(\alpha_1,\alpha_2,\ldots,\alpha_n)$, $\beta=(\beta_1,\beta_2,\ldots,\beta_n)\in\C^n$ such that $\alpha_j-\beta_k\notin\Z$ for any $1\le j,k\le n$, is defined as a subgroup  of $\GL_n(\C)$ generated by the companion matrices $A$ and $B$ of the polynomials 
\[f(x)=\prod_{j=1}^n(x-e^{2\pi i\alpha_j})\quad\mbox{ and }\quad g(x)=\prod_{j=1}^n(x-e^{2\pi i\beta_j}).\]

Levelt (cf. \cite[Theorem 3.5]{BH}) showed that there exists a basis of the (local) solution space of the hypergeometric differential equation
\[D(\alpha,\beta)u=0\]
defined on $\mathbb{P}^1\setminus\{0,1,\infty\}$, where $D(\alpha,\beta):=(\theta+\beta_1-1)\cdots(\theta+\beta_n-1)-z(\theta+\alpha_1)\cdots(\theta+\alpha_n)$
and $\theta=z\frac{d}{dz}$,  with respect to which the monodromy group of the hypergeometric equation is the hypergeometric group $\Gamma(\alpha,\beta)$ defined as above.

Now, we consider the hypergeometric groups $\Gamma(\alpha,\beta)$ for which the associated polynomials $f,g$ are products of cyclotomic polynomials and form a primitive pair (that is, both of the polynomials are, simultaneously, not polynomials in $x^k$ for any $k\ge 2$). Beukers and Heckman \cite[Theorem 6.5]{BH} show that in case when the constant terms of both polynomials are $1$ (this happens only when $n$, the degree of the polynomials $f,g$, is an even number), the corresponding hypergeometric group $\Gamma(\alpha,\beta)$ preserves a non-degenerate symplectic form $\Omega$ and $\Gamma(\alpha,\beta)$ is contained inside the integral symplectic group $\Sp_\Omega(\Z)$ as a Zariski dense subgroup. 

If we denote the Zariski closure of the hypergeometric group $\Gamma(\alpha,\beta)$ inside $\GL_n(\C)$ by $\mathbf{G}$ and $\Gamma(\alpha,\beta)\subseteq \mathbf{G}(\Z)$, then we call the hypergeometric group $\Gamma(\alpha,\beta)$ {\it arithmetic} if the index of $\Gamma(\alpha,\beta)$ inside $\mathbf{G}(\Z)$ is finite; and {\it thin}, otherwise.

The question of Sarnak \cite{S}, to determine the pairs of parameters $\alpha,\beta\in\Q^n$ for which the associated hypergeometric groups $\Gamma(\alpha,\beta)$ are arithmetic or thin, has been considered by many mathematicians and some progress on answering this question can be seen in \cite{BDSS, BS, BT, DFH, F, FMS, HvS, SS, SS0, SS17, SV, Ve2, Ve3}. The case when $\alpha=(0,0,0,0)$ has drawn much interest. In this case, there are 14 symplectic hypergeometric groups with a maximally unipotent monodromy associated to the pairs of polynomials $f, g$, where  $f(x)=(x-1)^4$ and $g(x)$ are products of cyclotomic polynomials such that $g(0)=1$, $g(1)\neq 0$ and the pair $f, g$ forms a primitive pair. These 14 hypergeometric groups arise as monodromy groups of families of Calabi-Yau threefolds fibering over $\mathbb{P}^1\setminus\{0,1,\infty\}$. It is now well known that half of these $14$ hypergeometric groups are arithmetic (cf.  \cite{SS} and \cite{SV}) and the other half are thin (cf. \cite{BT}). 

When $n\geq 6$ is an even number, consider $f(x)=(x-1)^n$ and $g(x)$ as described above. Then, one may ask if the hypergeometric groups associated with these pairs of polynomials also have the same dichotomy (as having half of them arithmetic and the other half thin) between arithmeticity and thinness that the $14$ hypergeometric groups (with a maximally unipotent monodromy) of degree four have. The first case to consider is $n=6$, and in this case a computation shows that there are 40 symplectic hypergeometric groups (cf. \cite[Table A]{BDSS}) associated with the pairs of polynomials $f(x)=(x-1)^6$ and $g(x)$ as described above. In \cite{BDSS}, we show that $18$ of these $40$ symplectic hypergeometric groups are arithmetic. 

\begin{rmk} It may be noted that according to \cite[Table 1]{BDN} the above dichotomy in $n=6$ is now disproved.
\end{rmk}

In this article, we show the thinness of $7$ of the above $40$ hypergeometric groups. 

\begin{thm}\label{maintheorem}
	The hypergeometric groups $\Gamma(\alpha,\beta)$ associated to the $7$ pairs of the parameters $\alpha,\beta$ where $\alpha=(0,0,0,0,0,0)$ and $\beta$ is any of the $7$ parameters appearing in Table \ref{Table-1} are thin.
\end{thm} 
\begin{table}[ht]
	\caption{Examples of thin hypergeometric groups in $\Sp(6)$}\label{Table-1}
	\centering
	\begin{tabular}{|c| c| c| c| c| c|}
	\hline
    S.No. & $\beta$ & $a$ & $d$ & $c$ & Splitting of $\Gamma(\alpha,\beta)$ \\ [0.5ex]
	\hline
	1 & (1/2,1/2,1/2,1/2,1/2,1/2) & 64 & 48 & 12 & $\mathbb{Z}*\mathbb{Z}$\\
	\hline
	2 & (1/2,1/2,1/2,1/2,1/3,2/3) & 48 & 40 & 11 & $\mathbb{Z}*\mathbb{Z}$\\
	\hline
	3 & (1/2,1/2,1/2,1/2,1/4,3/4) & 32 & 32 & 10 & $\mathbb{Z}*\mathbb{Z}$\\
	\hline
	4 & (1/2,1/2,1/2,1/2,1/6,5/6) & 16 & 29 & 9 & $\mathbb{Z}*\mathbb{Z}$\\
	\hline
	5 & (1/2,1/2,1/3,2/3,1/3,2/3) & 36 & 33 & 10 & $\mathbb{Z}*\mathbb{Z}$\\
	\hline
	6 & (1/2,1/2,1/3,2/3,1/4,3/4) & 24 & 26 & 9 & $\mathbb{Z}*\mathbb{Z}$\\
	\hline
	7 & (1/2,1/2,1/5,2/5,3/5,4/5) & 20 & 25 & 9 & $\mathbb{Z}*\mathbb{Z}$\\
	\hline
	\end{tabular}
	\vspace{0.5 cm}
	
	(The integers $a,d$ and $c$ determine the entries of some conjugate of the matrix $A$ (cf. Section \ref{notation}).)	
	\end{table}

We now sketch the proof of the above theorem. It follows from \cite[Table 2]{BCRR} and the universal coefficient theorem  for cohomology that $H^2(\Sp_6(\Z), \Q)$ is isomorphic to $\Q$ and hence the cohomological dimension of $\Sp_6(\Z)$ over $\Q$ is $\geq 2$. Since $\Sp_6(\Z)$ has no $\Q$-torsion, that  is, the order of every finite subgroup of $\Sp_6(\Z)$ is a unit in $\Q$, it follows from \cite[Theorem 9.2]{SW} that if $\Gamma$ is a finite index subgroup of $\Sp_6(\Z)$ then the cohomological dimension of $\Gamma$ over $\Q$ is also $\ge 2$. Therefore, to prove Theorem \ref{maintheorem} it is enough to prove that the hypergeometric groups appearing in Theorem \ref{maintheorem} are free groups or contain a free subgroup of finite index (free groups have cohomological dimensions $1$). 

Using the technique of Brav and Thomas \cite{BT}, we use the following version of the ping-pong lemma from \cite[Proposition III.12.4]{LS} to show that the hypergeometric groups $\Gamma(\alpha,\beta)$ appearing in Theorem \ref{maintheorem} are isomorphic to $\Z*\Z$.
\begin{thm}[Ping-Pong Lemma]\label{ping-pong lemma}
Let a group $G$ be generated by two of its subgroups, $G_1$ and $G_2$. Let $e$ be the identity element of the group $G$. Suppose at least one of these two subgroups has order greater than $2$ and $G$ acts on a set $W$. Suppose there are non-empty subsets $X$ and $Y$ of $W$ such that 
\begin{enumerate}
	\item[(i)] $X$ and $Y$ are disjoint subsets;
	\item[(ii)] $(G_2\setminus\{e\})X\subseteq Y$;
	\item[(iii)] $(G_1\setminus\{e\})Y \subseteq X$;
\end{enumerate}
Then $G=G_1* G_2$, the free product of $G_1$ and $G_2$.
\end{thm}

 For the hypergeometric groups appearing in Table \ref{Table-1}, we consider the standard action of $\Gamma(\alpha,\beta)$ on $\R^6$. We apply some change of basis in $\R^6$ to make computations simpler. Let $K$ be the change of the basis matrix.  Let $U=K^{-1}AK$ and $T=K^{-1}A^{-1}BK$, where $A$ and $B$ are the matrices defined at the beginning of this section. Let $R=TU$.

The hypergeometric group $\Gamma(\alpha,\beta)$ (with respect to the new basis) is generated by the cyclic subgroups $G_1=\langle T\rangle$ and $G_2=\langle R\rangle$. Since $A^{-1}B$ is a non-trivial unipotent matrix and $B$ has repeated roots in all the $7$ cases of Table \ref{Table-1},  $G_1$ and $G_2$ both are isomorphic to $\Z$.

We define the two sets $X$ and $Y$ as 
\[
X=\pm C^+\cup \pm C^-\qquad \mbox{ and } \qquad Y=\bigcup \limits_{j\in\Z\setminus\{0\}} R^jX,\]
where $C^+$ and $C^-$ are some cones in $\R^6$. With this setting we are able to verify the ping-pong conditions (i), (ii) and (iii) of Theorem \ref{ping-pong lemma} for the 7 cases of Table \ref{Table-1} and conclude that the corresponding hypergeometric groups $\Gamma(\alpha,\beta)$ are isomorphic to $\Z*\Z$.

\begin{rmk}
  Examples 1-7 appearing in Table \ref{Table-1} are, respectively, Examples 1, 2, 3, 4, 5, 6, 10 of \cite[Table A]{BDSS}. While we were trying to get all the thin hypergeometric groups of \cite[Table A]{BDSS} one by one, Bajpai, Dona and Nitsche announced their article \cite{BDN} on arXiv in which they also show the thinness of these examples. In comparison to \cite{BDN}, the added value of the present article is:
  \begin{enumerate}
  \item  It demonstrates that the specific method of Brav and Thomas works in higher dimensions with almost no changes.
		\item\label{remark(2)} The brute force  input to make the method work consists of only 2 real numbers appearing in $v$ (cf. Section \ref{notation}). Consequently, there is a chance that one can find a pattern in the choices of those numbers that generalizes to other hypergeometric groups.
  \end{enumerate}
\end{rmk}

\section{Notation and Strategy}\label{notation}
The hypergeometric group $\Gamma(\alpha,\beta)$ with the pair of parameters $\alpha=(0,0,0,0,0,0)$ and $\beta$ as in Table \ref{Table-1}, with respect to the new basis (obtained by using a change of basis matrix $K$, for suitable $K$, provided in Section \ref{Proof of freeness}), is generated by the matrices $U$, $T$ and $R=TU$, where $U=K^{-1}AK$ and $T=K^{-1}A^{-1}BK$. The matrices $U, T$ and $R$ have the following form and they preserve a symplectic form $J$ whose matrix form is as given below:

$$  U= 
\left[ \begin {array}{cccccc} 1&1&0&0&0&0\\ \noalign{\medskip}0&1&1&0
&0&0\\ \noalign{\medskip}0&0&1&0&0&0\\ \noalign{\medskip}-a&-a&0&1&0&0
\\ \noalign{\medskip}0&d&d&-1&1&0\\ \noalign{\medskip}0&0&-c&1&-1&1
\end {array} \right],\qquad T=\left[ \begin {array}{cccccc} 1&0&0&0&0&0\\ \noalign{\medskip}0&1&0&0
&0&0\\ \noalign{\medskip}0&0&1&0&0&1\\ \noalign{\medskip}0&0&0&1&0&0
\\ \noalign{\medskip}0&0&0&0&1&0\\ \noalign{\medskip}0&0&0&0&0&1
\end {array} \right],
$$

$$
R=\left[ \begin {array}{cccccc} 1&1&0&0&0&0\\ \noalign{\medskip}0&1&1&0
&0&0\\ \noalign{\medskip}0&0&1-c&1&-1&1\\ \noalign{\medskip}-a&-a&0&1&0
&0\\ \noalign{\medskip}0&d&d&-1&1&0\\ \noalign{\medskip}0&0&-c&1&-1&1
\end {array} \right]
,\qquad J=\left[ \begin {array}{cccccc} 0&0&0&1&0&0\\ \noalign{\medskip}0&0&0&0
&1&0\\ \noalign{\medskip}0&0&0&0&0&1\\ \noalign{\medskip}-1&0&0&0&0&0
\\ \noalign{\medskip}0&-1&0&0&0&0\\ \noalign{\medskip}0&0&-1&0&0&0
\end {array} \right].
$$

Observe that in case of the hypergeometric groups appearing in Theorem \ref{maintheorem}, $R^m\neq I$ for any $m\in\Z\setminus\{0\}$ as the characteristic polynomial of $B$ has repeated roots.

Following the idea of Brav and Thomas \cite{BT}, we consider first the matrix
$$ 
H=\left[ \begin {array}{cccccc} 1&1&0&0&0&0\\ \noalign{\medskip}0&-1&0&0
&0&0\\ \noalign{\medskip}0&0&1&0&0&-1\\ \noalign{\medskip}0&0&0&-1&0&0
\\ \noalign{\medskip}0&-d&0&-1&1&0\\ \noalign{\medskip}0&0&0&0&0&-1
\end {array} \right]
$$
where $d$ is defined as the $(5,2)$ entry of the matrix $U$. It follows that $H$ fixes the subspace $V$ spanned by its first, third and fifth column vectors. A computation shows that
\begin{equation}\label{H^2=I}
 H^2=I, \qquad HRH=R^{-1} \qquad \mbox{ and } \qquad HT^{-1}H=T
\end{equation}
where $I$ is the $6\times 6$ identity matrix.

For a $6\times 6$ unipotent matrix $F$, we define 
\[\log(F)=(F-I)-\frac{1}{2}(F-I)^2+\frac{1}{3}(F-I)^3-\frac{1}{4}(F-I)^4+\frac{1}{5}(F-I)^5.\]
Observe that $T^{-1}R$ and $TR^{-1}$ are unipotent matrices and define $P=\log(T^{-1}R)$ and $Q=\log(TR^{-1})$.

It follows from Equation (\ref{H^2=I}) that
\begin{equation*}
H(T^{-1}R)^j=(TR^{-1})^jH\qquad \mbox{ for all } j\in\N
\end{equation*}
and from this it follows that $HP=QH$, and we get 
\begin{equation}\label{HP^j=Q^jH}
 HP^j=Q^jH\qquad \mbox{ for all } j\in\N.
\end{equation}

Now, we choose a vector $v$ in $V$, and define 
\[C^+=\left\{\sum_{j=0}^5 t_jP^j v :  t_j\in\R_{> 0},\ \forall\ 0\le j\le 5\right\}\]the open cone generated by the vectors $v, Pv, P^2v, P^3v, P^4v, P^5v$, and
\[C^-=\left\{\sum_{j=0}^5 t_jQ^j v :  t_j\in\R_{> 0},\ \forall\ 0\le j\le 5\right\}\]the open cone generated by the vectors $v, Qv, Q^2v, Q^3v, Q^4v, Q^5v$ inside  $\R^6$. Since $HP^j=Q^jH$ (cf. Equation (\ref{HP^j=Q^jH})) and $Hv=v$, it follows that 
\begin{equation}
 C^-=HC^+\ .
\end{equation}

Now, we are ready to define the subsets $X$ and $Y$ of $\R^6$ (with the standard action of the hypergeometric group $\Gamma(\alpha,\beta)$) that satisfy the hypotheses of the ping-pong lemma (cf. Theorem \ref{ping-pong lemma}). Let
\begin{equation} 
X=\pm C^+\cup \pm C^-\qquad \mbox{ and } \qquad Y=\bigcup\limits_{j\in\Z\setminus\{0\}} R^jX\ .
\end{equation}

 Let  $M$ be the matrix whose column vectors are $v, Pv, P^2v, P^3v, P^4v, P^5v$, and let $N$ be the matrix whose column vectors are $v, Qv, Q^2v, Q^3v, Q^4v, Q^5v$. Observe that $N=HM$.

 \begin{rmk} For the vector $v$ we choose, the matrices $M$ and $N$ both are invertible and the choice of $v$ is motivated by the analog calculation
of \cite{BT} in which $v$ is the solution vector of the equation $v^tJPv={\bf{0}}$ in $V$. In \cite{BT} the vector $v$ is uniquely determined but in our case the equation $v^tJPv={\bf{0}}$ has infinitely many solution vectors in $V$ and not all solution vectors help in verifying the ping-pong conditions. Based on computations we choose, for each case in Table \ref{Table-1} separately, a suitable vector $v$ (constructed from the solution vectors of the equation $v^tJPv={\bf{0}}$).
 \end{rmk}

The hypergeometric group $\Gamma(\alpha,\beta)$ is generated by its subgroups $G_1=\langle T\rangle$ and $G_2=\langle R\rangle$. We show that $G_1\cap G_2=\{I\}$. Observe that $R^m=T^n$, for some $m,n\in\Z$, implies that $A^{-1}B^mA=C^n$, where $C=A^{-1}B$ fixes the standard basis vectors $e_1,e_2,\ldots,e_5$ in $\R^6$, and we get $A^{-1}B^mAe_j=C^ne_j=e_j$ for all $1\le j\le 5$.  This implies that $B^me_j=e_j$ for all $2\le j\le 6$. Also, $Be_1=e_2=B^me_2=B^m Be_1$ implies that $B^me_1=e_1$. It now follows that $B^m=I$. Since the characteristic polynomials of $B$ associated to the parameters appearing in Theorem \ref{maintheorem} have repeated roots, we get that $m=0=n$. Hence $G_1\cap G_2=\{I\}$.

Now, we need to verify the following conditions of the ping-pong lemma.
\begin{enumerate}
	\item[(A)] $X$ and $Y$ are disjoint subsets.
	\item[(B)] $(G_2\setminus\{I\})X\subseteq Y$.
	\item[(C)] $(G_1\setminus\{I\})Y \subseteq X$.
\end{enumerate}
We will verify condition $(A)$ by verifying the following statements.
\begin{enumerate}
\item[(A1)] $R^jC^+$ and $R^jC^-$ are disjoint from $\pm C^+$ for all $j \in \mathbb{Z}\setminus\{0\}$.
\item[(A2)] $R^jC^+$ and $R^jC^-$ are also disjoint from $\pm C^-$ for all $j \in \mathbb{Z}\setminus\{0\}$.
\end{enumerate}
Condition (B) follows from the construction of the subset $Y$. We will verify  condition $(C)$ by verifying the following statements.
\begin{enumerate}
\item[(C1)] $T^{-1}C^+ \subseteq C^+$.
\item[(C2)] $TC^- \subseteq C^-$.
\item[(C3)] $T^{-1}R^j(C^+ \cup C^-) \subseteq \pm C^+$ for all $j \in \mathbb{Z}\setminus\{0\}$.
\item[(C4)] $TR^j(C^+ \cup C^-) \subseteq \pm C^-$ for all $j \in \mathbb{Z}\setminus\{0\}$.
\end{enumerate}

Now, we have the following convention.
\begin{defn}
We call a row vector $(a_{i1}, a_{i2},\ldots, a_{in})$ of an $n\times n$ matrix $(a_{ij})$ non-negative (respectively, non-positive) if $a_{ij}\geq 0$ (respectively, $a_{ij}\leq 0$) for all $1\le j\le n$. Also, we call an $n\times n$ matrix $(a_{ij})$ non-negative (respectively, non-positive) if $a_{ij}\geq 0$  (respectively, $a_{ij}\leq 0$) for all $1\le i, j\le n$.
\end{defn}

Observe that $R^jC^+\cap C^+$ is non-empty exactly if there exist column vectors $v_1,v_2$ having all of their coordinates positive such that $R^jMv_1=Mv_2$, that is, $M^{-1}R^jMv_1=v_2$. Similarly, $R^jC^-\cap C^+$ is non-empty exactly if there exist column vectors $v_1,v_2$ having all of their coordinates positive  such that $R^jNv_1=Mv_2$, that is, $M^{-1}R^jNv_1=v_2$. Therefore, to verify the statement (A1), it is sufficient to show that the matrices $M^{-1}R^jM$ and $M^{-1}R^jN$, for all $j \in \mathbb{Z}\setminus\{0\}$, have some non-negative and some non-positive rows (as when this is the case, $M^{-1}R^jM$ and $M^{-1}R^jN$ map a column vector $v_1$ having all of its coordinates positive to a column vector $v_2$ having some positive and some negative coordinates). 
 
Since $N=HM$, $H^2=I$ and $HRH=R^{-1}$ (cf. Equation (\ref{H^2=I})), we get
\begin{equation}\label{2.5}
N^{-1}R^jM=(HM)^{-1}R^j(HN)=M^{-1}H^{-1}R^jHN=M^{-1}HR^jHN=M^{-1}R^{-j}N
 \end{equation}
 \begin{equation}\label{2.6} N^{-1}R^jN=(HM)^{-1}R^j(HM)=M^{-1}H^{-1}R^jHM=M^{-1}HR^jHM=M^{-1}R^{-j}M.
   \end{equation}

It follows from Equations (\ref{2.5}) and (\ref{2.6}) that if the matrices $M^{-1}R^jM$ and $M^{-1}R^jN$, for all $j \in \mathbb{Z}\setminus\{0\}$, have some non-negative and some non-positive rows, then the matrices $N^{-1}R^jN$ and $N^{-1}R^jM$, for all $j \in \mathbb{Z}\setminus\{0\}$, also have some non-negative and some non-positive rows.  This shows that, for all $j\in\Z\setminus\{0\}$, $R^jC^+$ and $R^jC^-$ both are disjoint from $\pm C^-$. Hence the sufficient criterion for (A1) above also proves (A2).
 
The condition (C) is divided into four parts from (C1) to (C4). Observe that $T^{-1}C^+\subseteq C^+$ if for any column vector $v_1$ having all of its coordinates some positive real numbers, there exists a column vector $v_2$ having all of its coordinates some positive real numbers such that $T^{-1}Mv_1=Mv_2$. Thus, (C1) is equivalent to showing that the matrix $M^{-1}T^{-1}M$ is non-negative (as when this is the case, $M^{-1}T^{-1}M$ will map a column vector $v_1$ having all of its coordinates positive to a column vector $v_2$ having all of its coordinates positive).

Similarly, $TC^-\subseteq C^-$ if for any column vector $v_1$ having all of its coordinates some positive real numbers, there exists a column vector $v_2$ having all of its coordinates some positive real numbers such that $TNv_1=Nv_2$. Thus, (C2) is equivalent to showing that the matrix $N^{-1}TN$ is non-negative. Again, since $N=HM$, $H^2=I$ and $HTH=T^{-1}$ (cf. Equation (\ref{H^2=I})), we get
\begin{equation}N^{-1}TN=(HM)^{-1}T(HM)=M^{-1}H^{-1}THM=M^{-1}(HTH)M=M^{-1}T^{-1}M\end{equation}
and hence the sufficient criterion for (C1) also proves (C2).

By using similar arguments as above, to verify the statement (C3) it is enough to show that both the matrices $M^{-1}T^{-1}R^jM$ and $M^{-1}T^{-1}R^jN$, for all $j \in \mathbb{Z}\setminus\{0\}$, are either non-negative or non-positive. 

To show (C4), it is enough to show that both the matrices $N^{-1}TR^jM$ and $N^{-1}TR^jN$, for all $j \in \mathbb{Z}\setminus\{0\}$, are either non-negative or non-positive. Again, since $N=HM$, $H^2=I$, $HRH=R^{-1}$ and $HTH=T^{-1}$ (cf. Equation (\ref{H^2=I})),  we get  		
\begin{equation}
N^{-1}TR^jM=(HM)^{-1}TR^j(HN)=M^{-1}H^{-1}TR^jHN=M^{-1}HTR^jHN=M^{-1}T^{-1}R^{-j}N
\end{equation} 		
\begin{equation}
N^{-1}TR^jN=(HM)^{-1}TR^j(HM)=M^{-1}HTR^jHM=M^{-1}T^{-1}R^{-j}M
\end{equation}
and hence the sufficient criterion for (C3) also proves (C4).
 		
 
 In conclusion, to apply the ping-pong lemma to the groups in Table \ref{Table-1} with $X$ and $Y$ obtained as above from a given $v\in\R^6$, it suffices to check that $M$ has full rank, that $M^{-1}R^jM$ and $M^{-1}R^jN$ have both non-positive and non-negative rows for all $j\neq 0$ (A1), that $M^{-1}T^{-1}M$ is non-negative (C1) and that $M^{-1}T^{-1}R^jM$ and $M^{-1}T^{-1}R^jN$ are non-negative or non-positive for all $j\neq 0$ (C3).
 
\section{Proof of the freeness of the hypergeometric groups of Table \ref{Table-1}}\label{Proof of freeness}

As it is explained above in Section \ref{notation}, to show that the hypergeometric groups associated with the parameters appearing in Table \ref{Table-1} are free groups (using the ping-pong lemma), we need to verify only the statements (A1), (C1) and (C3) of Section \ref{notation} for these groups. 

We write below the detailed verification for Example 1, and also for Example 2 as it requires some different explanations. For other examples, we provide the needed data which helps in verifying the ping-pong conditions using the methods of Examples 1 and 2. The computations are an adaptation of the methods
in \cite{BT} to our setting. The main new ingredient is the choice of v. See Section \ref{notation} for the notations.


\subsection{Case 1}\label{case1} \textbf{$\alpha=(0,0,0,0,0,0), \quad \beta=(1/2,1/2,1/2,1/2,1/2,1/2)$}.

In this case the corresponding polynomials $f,g$ are 
\[f(x)={x}^{6}-6\,{x}^{5}+15\,{x}^{4}-20\,{x}^{3}+15\,{x}^{2}-6\,x+1\quad\mbox{ and }\quad g(x)={x}^{6}+6\,{x}^{5}+15\,{x}^{4}+20\,{x}^{3}+15\,{x}^{2}+6\,x+1\] 
and the corresponding hypergeometric group $\Gamma(\alpha,\beta)$ is generated by the companion matrices $A, B$ of the polynomials $f,g$. We first consider the following change of basis matrix
 \[K=\left[ \begin {array}{cccccc} 0&48&-12&0&-1&-1\\ \noalign{\medskip}-
64&-144&0&-1&4&5\\ \noalign{\medskip}128&208&-40&3&-6&-10
\\ \noalign{\medskip}-64&-112&0&-3&4&10\\ \noalign{\medskip}0&0&-12&1&
-1&-5\\ \noalign{\medskip}0&0&0&0&0&1\end {array} \right].\]

Then $U=K^{-1}AK$, $T=K^{-1}A^{-1}BK$ and $R=TU$ have the form described in Section 2 with  $a=64$, $d=48$ and $c=12$.

Guided by computer calculations we pick  $v=\left(-\frac{1}{10}, 0,1, 0, -\frac{1162}{225}, 0\right)$ inside the fixed subspace $V$ of the matrix $H$. Note that in this case $-R$ is unipotent and so we can define $Z=\log(-R)$. Then, as required, the resulting matrix $M$ has full rank, and
\[-R= e^Z= I+Z+\frac{Z^2}{2!}+\frac{Z^3}{3!}+\frac{Z^4}{4!}+\frac{Z^5}{5!}\]
\begin{equation}\label{Rj} R^j= (-1)^je^{jZ}= (-1)^j\left(I+jZ+\frac{j^2Z^2}{2!}+\frac{j^3Z^3}{3!}+\frac{j^4Z^4}{4!}+\frac{j^5Z^5}{5!}\right).
\end{equation}
To verify (A1) we compute the matrices $M^{-1}R^jM$ and $M^{-1}R^jN$, for all $j \in \mathbb{Z}\setminus\{0\}$, and show that they have some non-negative and some non-positive rows.

Denote $A_i=M^{-1}\frac{Z^i}{i!}M$, for $i=1,2,3,4,5$. Then
\begin{equation}\label{disjointXY_1}
M^{-1}R^jM= (-1)^j(I+jA_1+j^2A_2+j^3A_3+j^4A_4+j^5A_5)\ .
\end{equation}
By computation we get 
\begin{equation}\label{case1-A5}
A_5=\left[ \begin {array}{cccccc} {\frac{98}{3375}}&{\frac{916}{3375}}&{
	\frac{304}{225}}&{\frac{256}{75}}&{\frac{64}{15}}&{\frac{128}{15}}
\\ \noalign{\medskip}-{\frac{49}{3375}}&-{\frac{458}{3375}}&-{\frac{
		152}{225}}&-{\frac{128}{75}}&-{\frac{32}{15}}&-{\frac{64}{15}}
\\ \noalign{\medskip}{\frac{49}{16875}}&{\frac{458}{16875}}&{\frac{152
	}{1125}}&{\frac{128}{375}}&{\frac{32}{75}}&{\frac{64}{75}}
\\ \noalign{\medskip}-{\frac{49}{202500}}&-{\frac{229}{101250}}&-{
	\frac{38}{3375}}&-{\frac{32}{1125}}&-{\frac{8}{225}}&-{\frac{16}{225}}
\\ \noalign{\medskip}{\frac{8477}{12150000}}&{\frac{39617}{6075000}}&{
	\frac{3287}{101250}}&{\frac{1384}{16875}}&{\frac{346}{3375}}&{\frac{
		692}{3375}}\\ \noalign{\medskip}-{\frac{8477}{24300000}}&-{\frac{39617
	}{12150000}}&-{\frac{3287}{202500}}&-{\frac{692}{16875}}&-{\frac{173}{
		3375}}&-{\frac{346}{3375}}\end {array} \right].
		\end{equation}

Notice that $A_5$ has some non-negative and some non-positive rows. Now, we use the Archimedean property of the real numbers to show that there exists a positive integer $n\in\Z$ such that the entries of the first two rows of the matrix $M^{-1}R^jM$, for all $j \in \mathbb{Z}\setminus\{0\}$ and $|j|\geq n$,  have the same sign as that of the corresponding entries of the matrix  $\mathrm{sign}(j)(-1)^{|j|}A_5$.

A sufficient $n$ can be computed after writing
 \begin{equation}	(-1)^j(I+jA_1+j^2A_2+j^3A_3+j^4A_4+j^5A_5)=(-1)^j(I+j(A_1+j(A_2+j(A_3+j(A_4+jA_5))))) .\end{equation}
 
By using the Archimedean property of the real numbers on the entries of the matrices $A_5$ and $A_4$, we get a positive integer $k_5$ such that the entries of the matrix $A_4+jA_5$ have the same sign as that of the corresponding entries of the matrix $A_5$ for all $j\geq k_5$. Again using the Archimedean property of the real numbers on the entries of the matrices $A_4+k_5A_5$ and $A_3$, we get a positive integer $k_4$ such that the entries of the matrix $A_3+j(A_4+k_5A_5)$ have the same sign as that of the corresponding entries of the matrix $A_4+k_5A_5$ for all $j\geq k_4$.  Repeating this process we get the positive integers $k_5,k_4,k_3,k_2,k_1$ such that for all $j\ge n=\max\{k_5,k_4,k_3,k_2,k_1\}$ the entries of the matrix $(-1)^j(1+j(A_1+j(A_2+j(A_3+j(A_4+jA_5)))))$ have the same sign as that of the corresponding entries of the matrix $(-1)^jA_5$.

In case when $j$ is a negative integer, we take $j=-l$ where $l$ is a positive integer and write the expression for the matrix $M^{-1}R^jM$ as
\begin{equation}
	(-1)^j(1+jA_1+j^2A_2+j^3A_3+j^4A_4+j^5A_5)=(-1)^l(1+l(-A_1+l(A_2+l(-A_3+l(A_4-lA_5)))))
\end{equation}
and use the Archimedean property of the real numbers, as it is described above, to get a positive integer $m$ so that the matrix $M^{-1}R^jM$, for all $j\le -m  \in \mathbb{Z}\setminus\{0\}$, has some non-negative and some non-positive rows.
In this case, we get that $\max\{n,m\}$ is $5$. So, it follows that for $|j|\ge 5$, the matrix $M^{-1}R^jM$ has some non-negative and some non-positive rows.

For the remaining values of $j \in \mathbb{Z}\setminus\{0\}$, that is, for $j\in\Z$ with $1\leq |j|\leq 4$, we compute the matrices $M^{-1}R^jM$ individually and verify that the matrix $M^{-1}R^jM$, also for $1\leq |j|\leq 4$, has some non-negative and some non-positive rows.

Similarly, we verify that the matrix $M^{-1}R^jN$, for all $j \in \mathbb{Z}\setminus\{0\}$, has some non-negative and some non-positive rows. 

To verify (C1) we need to show that  the matrix $M^{-1}T^{-1}M$ is non-negative. By computation we get
\[M^{-1}T^{-1}M=\left[ \begin {array}{cccccc} 1&{\frac{458}{225}}&0&{\frac{128}{5}}&0
&64\\ \noalign{\medskip}0&1&0&0&0&0\\ \noalign{\medskip}0&{\frac{229}{
		1125}}&1&{\frac{64}{25}}&0&{\frac{32}{5}}\\ \noalign{\medskip}0&0&0&1&0
&0\\ \noalign{\medskip}0&{\frac{39617}{810000}}&0&{\frac{692}{1125}}&1
&{\frac{346}{225}}\\ \noalign{\medskip}0&0&0&0&0&1\end {array} \right].\]

To verify (C3) it is sufficient to show that both the matrices $M^{-1}T^{-1}R^jM$ and $M^{-1}T^{-1}R^jN$, for all $j \in \mathbb{Z}\setminus\{0\}$, are either non-negative or non-positive. 

Using Equation (\ref{Rj}) for the expression of $R^j$, we get
\begin{equation}
M^{-1}T^{-1}R^jM= (-1)^j(M^{-1}T^{-1}M+jD_1+j^2D_2+j^3D_3+j^4D_4+j^5D_5)
\end{equation}
\begin{equation}
M^{-1}T^{-1}R^jN= (-1)^j(M^{-1}T^{-1}N+jE_1+j^2E_2+j^3E_3+j^4E_4+j^5E_5)
\end{equation}
where $D_i=M^{-1}\frac{T^{-1}Z^i}{i!}M$ and $E_i=M^{-1}\frac{T^{-1}Z^i}{i!}N$, for $i=1,2,3,4,5$. 

By computation, we get 
\[D_5=E_5=\left[ \begin {array}{cccccc} -{\frac{98}{3375}}&-{\frac{916}{3375}}&
-{\frac{304}{225}}&-{\frac{256}{75}}&-{\frac{64}{15}}&-{\frac{128}{15}
}\\ \noalign{\medskip}-{\frac{49}{3375}}&-{\frac{458}{3375}}&-{\frac{
		152}{225}}&-{\frac{128}{75}}&-{\frac{32}{15}}&-{\frac{64}{15}}
\\ \noalign{\medskip}-{\frac{49}{16875}}&-{\frac{458}{16875}}&-{\frac{
		152}{1125}}&-{\frac{128}{375}}&-{\frac{32}{75}}&-{\frac{64}{75}}
\\ \noalign{\medskip}-{\frac{49}{202500}}&-{\frac{229}{101250}}&-{
	\frac{38}{3375}}&-{\frac{32}{1125}}&-{\frac{8}{225}}&-{\frac{16}{225}}
\\ \noalign{\medskip}-{\frac{8477}{12150000}}&-{\frac{39617}{6075000}}
&-{\frac{3287}{101250}}&-{\frac{1384}{16875}}&-{\frac{346}{3375}}&-{
	\frac{692}{3375}}\\ \noalign{\medskip}-{\frac{8477}{24300000}}&-{\frac
	{39617}{12150000}}&-{\frac{3287}{202500}}&-{\frac{692}{16875}}&-{\frac
	{173}{3375}}&-{\frac{346}{3375}}\end {array} \right].\]

Since all entries of the matrices $D_5$ and $E_5$ have the same signs, by using the Archimedean property of the real numbers, as above, we find that for $|j|\geq 7$ the matrices $M^{-1}T^{-1}R^jM$ and $M^{-1}T^{-1}R^jN$ are either non-negative or non-positive. We compute that the matrices $M^{-1}T^{-1}R^jM$ and $M^{-1}T^{-1}R^jN$, also for  $1\le |j|\le 6$,  are either non-negative or non-positive. 

This completes the proof of the freeness of the hypergeometric group.

\subsection{Case 2}\textbf{$\alpha=(0,0,0,0,0,0),\quad \beta=(1/2,1/2,1/2,1/2,1/3,2/3)$}.

We consider the change of basis matrix
\[K=\left[ \begin {array}{cccccc} 0&40&-11&0&-1&-1\\ \noalign{\medskip}-
48&-120&4&-1&4&5\\ \noalign{\medskip}96&168&-34&3&-6&-10
\\ \noalign{\medskip}-48&-88&4&-3&4&10\\ \noalign{\medskip}0&0&-11&1&-
1&-5\\ \noalign{\medskip}0&0&0&0&0&1\end {array} \right].\]

In this case $a=48$, $d=40$, $c=11$ and  $v=\left(-\frac{1}{10},0,1,0,-\frac{631}{150},0\right)$. This time neither $R$ nor $-R$ is unipotent. But $R^6$ is unipotent, so let $Z=\log(R^6)$. Note that $Z^4=0$. Then
\[R^6=e^Z=I+Z+\frac{Z^2}{2!}+\frac{Z^3}{3!},\]
\begin{equation}\label{R6n+k}
R^{6n+k}=R^kR^{6n}=R^ke^{nZ}=R^k\left(I+nZ+\frac{n^2Z^2}{2!}+\frac{n^3Z^3}{3!}\right).
\end{equation}

We write $j=6n+k$ for $n \in \mathbb{N}\cup\{0\}$ and $0\leq k \leq 5$. If we denote $A_{i,k}=M^{-1}\frac{R^kZ^i}{i!}M$  and $B_{i,k}=M^{-1}\frac{R^kZ^i}{i!}N$   for $i=1,2,3$, then
\begin{equation}{\label{R6n+k with M}}
	M^{-1}R^{6n+k}M=M^{-1}R^kM+nA_{1,k}+n^2A_{2,k}+n^3A_{3,k}
	\end{equation}
	\begin{equation}\label{R6n+k with N}
	M^{-1}R^{6n+k}N=M^{-1}R^kN+nB_{1,k}+n^2B_{2,k}+n^3B_{3,k}\ .
	\end{equation}
 
 By denoting $D_{i,k}=M^{-1}\frac{T^{-1}R^kZ^i}{i!}M$ and $E_{i,k}=M^{-1}\frac{T^{-1}R^kZ^i}{i!}N$ for $i=1,2,3$ and using the above Equation (\ref{R6n+k}), we get
 \begin{equation}{\label{6n+k with T^{-1} and M}}
 	M^{-1}T^{-1}R^{6n+k}M=M^{-1}T^{-1}R^kM+nD_{1,k}+n^2D_{2,k}+n^3D_{3,k}
 \end{equation}
\begin{equation}{\label{6n+k with T^{-1} and N}}
	M^{-1}T^{-1}R^{6n+k}N=M^{-1}T^{-1}R^kN+nE_{1,k}+n^2E_{2,k}+n^3E_{3,k}\ .
\end{equation}

For each $k$ we can find a $q\in\N$ such that the signs of the entries of $M^{-1}R^{6n+k}M$, $M^{-1}R^{6n+k}N$, $M^{-1}T^{-1}R^{6n+k}M$ and $M^{-1}T^{-1}R^{6n+k}N$ are the same as that of the entries of $A_{3,k}, B_{3,k}, C_{3,k}$ and $D_{3,k}$ for all $|n| > q$, and we verify the ping-pong conditions as in Case 1.

\subsection{Case 3}\textbf{$\alpha=(0,0,0,0,0,0),\quad \beta=(1/2,1/2,1/2,1/2,1/4,3/4)$}.


 We consider the change of basis matrix
 $$K=\scriptsize \left[ \begin {array}{cccccc} 0&32&-10&0&-1&-1\\ \noalign{\medskip}-
 32&-96&8&-1&4&5\\ \noalign{\medskip}64&128&-28&3&-6&-10
 \\ \noalign{\medskip}-32&-64&8&-3&4&10\\ \noalign{\medskip}0&0&-10&1&-
 1&-5\\ \noalign{\medskip}0&0&0&0&0&1\end {array} \right].
 $$
 
 
In this case $a=32$, $d=32$, $c=10$, $v=\left(-\frac{1}{11},0,1,0,-{\frac{277}{75}},0\right)$ and $Z=\log(R^4)$. 

\subsection{Case 4}\textbf{$\alpha=(0,0,0,0,0,0),\quad \beta=(1/2,1/2,1/2,1/2,1/6,5/6)$}.


 We consider the change of basis matrix
$$K= \scriptsize \left[ \begin {array}{cccccc} 0&24&-9&0&-1&-1\\ \noalign{\medskip}-16
&-72&12&-1&4&5\\ \noalign{\medskip}32&88&-22&3&-6&-10
\\ \noalign{\medskip}-16&-40&12&-3&4&10\\ \noalign{\medskip}0&0&-9&1&-
1&-5\\ \noalign{\medskip}0&0&0&0&0&1\end {array} \right].
$$


In this case $a=16$, $d=24$, $c=9$, $v=\left(-\frac{1}{10},0,1,0,-{\frac{3041}{810}},0\right)$ and $Z=\log(-R^3)$. 


\subsection{Case 5}\textbf{$\alpha=(0,0,0,0,0,0),\quad \beta=(1/2,1/2,1/3,2/3,1/3,2/3)$}.


 We consider the change of basis matrix
$$K=\scriptsize \left[ \begin {array}{cccccc} 0&33&-10&0&-1&-1\\ \noalign{\medskip}-
36&-99&7&-1&4&5\\ \noalign{\medskip}72&135&-30&3&-6&-10
\\ \noalign{\medskip}-36&-69&7&-3&4&10\\ \noalign{\medskip}0&0&-10&1&-
1&-5\\ \noalign{\medskip}0&0&0&0&0&1\end {array} \right].
$$

In this case $a=36$, $d=33$, $c=10$, $v=\left(-\frac{1}{6},0,1,0,-{\frac{207}{40}},0\right)$ and $Z=\log(R^6)$.

\subsection{Case 6}\textbf{$\alpha=(0,0,0,0,0,0),\quad \beta=(1/2,1/2,1/3,2/3,1/4,3/4)$}.


We consider the change of basis matrix
$$K= \scriptsize \left[ \begin {array}{cccccc} 0&26&-9&0&-1&-1\\ \noalign{\medskip}-24
&-78&10&-1&4&5\\ \noalign{\medskip}48&102&-26&3&-6&-10
\\ \noalign{\medskip}-24&-50&10&-3&4&10\\ \noalign{\medskip}0&0&-9&1&-
1&-5\\ \noalign{\medskip}0&0&0&0&0&1\end {array} \right].
$$

In this case $a=24$, $d=26$, $c=9$, $v=\left(-\frac{1}{6},0,1,0,-{\frac{472}{105}},0\right)$ and $Z=\log(R^{12})$. 




\subsection{Case 7}\textbf{$\alpha=(0,0,0,0,0,0),\quad \beta=(1/2,1/2,1/5,2/5,3/5,4/5)$}.


 We consider the change of basis matrix
 $$K= \scriptsize \left[ \begin {array}{cccccc} 0&25&-9&0&-1&-1\\ \noalign{\medskip}-20
 &-75&11&-1&4&5\\ \noalign{\medskip}40&95&-24&3&-6&-10
 \\ \noalign{\medskip}-20&-45&11&-3&4&10\\ \noalign{\medskip}0&0&-9&1&-
 1&-5\\ \noalign{\medskip}0&0&0&0&0&1\end {array} \right] .
 $$

 In this case $a=20$, $d=25$, $c=9$, $v=\left(-\frac{1}{10},0,1,0,-{\frac{221}{60}},0\right)$ and $Z=\log(R^{10})$.
 
 
 

 \section*{Acknowledgements}
We thank the referee for the invaluable comments and suggestions, which substantially improved the article's presentation. The work of the first named author is supported by the MATRICS grant no. MTR/2021/000368. This paper is part of the Ph.D. thesis of the second named author who would like to thank IIT Bombay for providing support. 

\section*{Competing Interests}
 The authors declare none.

\bibliography{Sp6}
\end{document}